\documentclass[graybox]{svmult}

\usepackage{type1cm}        % activate if the above 3 fonts are
\usepackage{makeidx}         % allows index generation
\usepackage{graphicx}        % standard LaTeX graphics tool
\usepackage{multicol}        % used for the two-column index
\usepackage[bottom]{footmisc}% places footnotes at page bottom

\usepackage{newtxtext}       %
\usepackage{newtxmath}       % selects Times Roman as basic font

\makeindex             % used for the subject index
\usepackage{booktabs}
\usepackage{url}

\newcommand{\bsc}{\boldsymbol{c}}
\newcommand{\bsk}{\boldsymbol{k}}
\newcommand{\bsu}{\boldsymbol{u}}
\newcommand{\bsx}{\boldsymbol{x}}

\newcommand{\rd}{\,\mathrm{d}}
\newcommand{\dunif}{\mathsf{U}}

\newcommand{\err}{\varepsilon}
\newcommand{\e}{\mathbb{E}}
\newcommand{\var}{\mathrm{var}}
\newcommand{\rmse}{\mathrm{RMSE}}
\begin{document}

\title*{On dropping the first Sobol' point}
\author{Art B. Owen}
\institute{Art B. Owen \at Stanford University, Stanford CA 94305, \email{owen@stanford.edu}
}
\maketitle

\abstract{
Quasi-Monte Carlo (QMC) points are a substitute for plain Monte Carlo (MC)
points that greatly improve integration accuracy under
mild assumptions on the problem.
Because QMC can give errors that are $o(1/n)$ as $n\to\infty$,
and randomized versions can attain root mean squared
errors that are $o(1/n)$,
changing even one point can change the estimate by an
amount much larger than the error would have been
and worsen the convergence rate.
As a result, certain practices that fit quite naturally
and intuitively with MC
points can be very detrimental to QMC performance.
These include thinning, burn-in, and taking sample
sizes such as powers of $10$, when
the QMC points were designed for different sample sizes.
This article looks at the effects of a common practice
in which one skips the first point of a Sobol' sequence.
The retained points ordinarily fail to be a digital net and
when scrambling is applied, skipping over the first
point can increase the numerical error by a factor
proportional to $\sqrt{n}$ where $n$ is the
number of function evaluations used.
}

\section{Introduction}
A Sobol' sequence
is an infinite sequence of points $\bsu_1,\bsu_2,\dots\in[0,1]^d$
constructed to fill out the unit cube with low
discrepancy, meaning that a measure of the distance between the
discrete uniform distribution on $\bsu_1,\dots,\bsu_n$
and the continuous uniform distribution on $[0,1]^d$ is made
small.
These points are ordinarily used to approximate
$$
\mu = \int_{[0,1]^d}f(\bsx)\rd\bsx
\quad\text{by}\quad
\hat\mu = \hat\mu_{\bsu,1} = \frac1n\sum_{i=1}^n
f(\bsu_i).
$$
The reason for calling this estimate $\hat\mu_{\bsu,1}$ will become apparent later.
Sobol' sequences are often used to estimate expectations with
respect to unbounded random variables, such as Gaussians.
In such cases $f$ subsumes a transformation from the uniform
distribution on $[0,1]^d$ to some other more appropriate distribution.
This article uses $1$-based indexing, so that the initial point is
$\bsu_1$. Sometimes $0$-based indexing is used, and
then the initial point is denoted $\bsu_0$.  Both
indexing conventions are widespread in mathematics and software
for Sobol' points and both have their benefits.  Whichever convention
is used, the first point should not be dropped.

%Taking $n=2^m$ for an integer $m\geqslant0$ is a best practice because then
%the points $\bsu_1,\dots,\bsu_n$ comprise a digital net as
%defined below.

The initial point of the Sobol' sequence is $\bsu_1=(0,0,\dots,0)$.
A common practice is to skip that point, similar to the burn-in
practice in Markov chain Monte Carlo (MCMC).
One then estimates $\mu$ by
$$
\hat\mu = \hat\mu_{\bsu,2} = \frac1n\sum_{i=2}^{n+1} f(\bsu_i).
$$
One reason to skip the initial point is that a transformation to a
Gaussian distribution might make the initial Gaussian point infinite.
That is problematic not just for integration problems
but also when $f$ is to be evaluated
at the design points to create surrogate models
for Bayesian optimization \cite{fraz:2018,bala:etal:2019}.
If one skips the initial point, then the next point in a Sobol'
sequence is usually $(1/2,1/2,\dots,1/2)$. While that is
an intuitively much more reasonable place to start,
starting there has detrimental consequences and there are better
remedies, described here.

A discussion about whether to  drop the initial point came up in the plenary
tutorial of Fred Hickernell at MCQMC 2020
about QMCPy \cite{QMCPy} software for QMC,
discussed in \cite{ChoEtal22a}.
%\url{https://github.com/QMCSoftware/QMCSoftware}.
The issue has been
discussed by the pytorch \cite{paszke2019pytorch}
community at
\url{https://github.com/pytorch/pytorch/issues/32047},
and the scipy \cite{virtanen2020scipy}  community   at
\url{https://github.com/scipy/scipy/pull/10844},
which are both incorporating QMC methods.
QMC and RQMC code for scipy is documented
at \url{https://scipy.github.io/devdocs/reference/stats.qmc.html}.

This article shows that skipping even one point of the
Sobol' sequence can be very detrimental.
The resulting points are no longer a digital net
in general, and in the case of scrambled Sobol'
points, skipping a point can bring about an inferior rate of convergence,
making the estimate less accurate by a factor that is
roughly proportional to~$\sqrt{n}$.

A second difficulty with Sobol' sequence
points is that  it is difficult to estimate the size  $|\hat\mu-\mu|$ of the integration
error from the data.
The well-known Koksma-Hlawka
inequality \cite{hick:2014} bounds $|\hat\mu-\mu|$ by
the product of two unknown quantities that are
extremely hard to compute, and while tight for some worst case integrands, it can
yield an extreme overestimate of the error,
growing ever more conservative as the dimension $d$ increases.

Randomly scrambling the Sobol' sequence points
preserves their balance properties and provides a basis
for uncertainty quantification.
Scrambling turns points $\bsu_i$ into random points $\bsx_i\sim \dunif[0,1]^d$.
The points $\bsx_1,\dots,\bsx_n$ are not independent.
Instead they retain the digital net property of Sobol' points
and consequent accuracy properties.
The result is randomized QMC (RQMC) points.
RQMC points also have some additional accuracy properties stemming
from the randomization.
With scrambled Sobol' points, we estimate $\mu$ by
$$
\hat\mu = \hat\mu_{\bsx,1} = \frac1n\sum_{i=1}^n f(\bsx_i).
$$
One can estimate the mean squared
error using $R$ independent replicates of
the $n$-point RQMC estimate $\hat\mu_{\bsx,1}$.
It is also possible to drop the first point in
RQMC, estimating $\mu$ by
\begin{align*}
\hat\mu=\hat\mu_{\bsx,2} = \frac1n\sum_{i=2}^{n+1} f(\bsx_i).
\end{align*}
The purpose of this paper is to show that $\hat\mu_{\bsx,1}$
is a much better choice than $\hat\mu_{\bsx,2}$.

Many implementations of a Sobol' sequence will produce
$n=2^m$ points $\bsu_i\in \{ 0, 1/n, 2/n,\dots,(n-1)/n\}^d\subset[0,1)^d$.
In that case, there is a safer way to avoid having a point at the origin
than skipping the first point.  We can use $\bsu_i+1/(2n)$ componentwise
and still have a digital net.
This is reasonable if we have already decided on the
value of $n$ to use.
It does not work to add that same value $1/(2n)$ to the next $2^m$ points and subsequent
values.  For one thing, the result may produce values
on the upper boundary of $[0,1]^d$ in the very next
batch and will eventually place points outside of $[0,1]^d$.
It remains better to scramble the Sobol' points.

We will judge the accuracy of integration via scrambled Sobol' points
through $\e( (\hat\mu-\mu)^2)^{1/2}$,  the root mean squared error (RMSE).
Plain Monte Carlo (MC) attains an RMSE of $O(n^{-1/2})$ for integrands $f\in L^2[0,1]^d$.

This paper is organized as follows.
Section~\ref{sec:nets} defines digital nets
and shows that skipping over the first point
can destroy the digital net property underlying
the analysis of Sobol' sequences. It also presents
properties of scrambled digital nets.
Section~\ref{sec:synth} shows some empirical
investigations on some very simple and favorable
integrands where  $\hat\mu_{\bsx,1}$
has an RMSE very near to the rate $O(n^{-3/2})$
while $\hat\mu_{\bsx,2}$ has an RMSE
very near to $O(n^{-1})$. These are both in line with
what we expect from asymptotic theory.
The relevance is not that our integrands are as trivial
as those examples, but rather that when realistic integrands
are well approximated by such simple ones
we get accuracy comparable to using those simple functions as control variates
\cite{hick:lemi:owen:2005}  but without
us having to search for control variates.
In skipping the first point we stand
to lose a lot of accuracy in integrating the simple functions and others close to them.
There is also no theoretical reason to expect
$\hat\mu_{\bsx,2}$ to have a smaller RMSE
than  $\hat\mu_{\bsx,1}$ does,
and so there is a Pascal's wager argument against
dropping the first point.
Section~\ref{sec:wingweight} looks at a ten dimensional
function representing the weight of an airplane wing
as a function of the way it was made. We see there
that skipping the first point is very detrimental.
Section~\ref{sec:discussion} considers some very special cases
where burn-in might be harmless,
recommends against using round number sample sizes
and thinning for QMC points, and discusses
distributing QMC points over multiple parallel processors.

\section{Digital nets and scrambling}\label{sec:nets}

In this section we review digital nets and describe properties
of their scrambled versions.
The points from Sobol' sequences provide the most
widely used example of digital nets.
For details of their construction and analysis,
see the monographs \cite{dick:pill:2010,nied:1992}.
There are numerous implementations
of Sobol' sequences \cite{brat:fox:1988,joe:kuo:2008,sobo:asot:krei:kuch:2011}.
They differ in what are called
`direction numbers' and they can also vary in the
order with which the points are generated. The numerical
results here use direction numbers
from \cite{joe:kuo:2008} via an implementation
from Nuyens' magic point shop, described in \cite{kuo:nuye:2016}
and scrambled as in \cite{rtms}.
The Sobol' and scrambled Sobol'  points in this paper were generated
using the R function {\tt rsobol} that appears in
\url{http://statweb.stanford.edu/~owen/code/}
along with some documentation.
That code also includes the faster and more space efficient
scrambling of Matousek \cite{mato:1998:2}.
%The generating matrices are from Dirk Nuyens' magic
%point shop following the description in \cite{kuo:nuye:2016}
%and based on direction numbers from \cite{joe:kuo:2008}.

We begin with the notion of elementary intervals,
which are special hyper-rectangular subsets of $[0,1)^d$.
For an integer base $b\geqslant2$, a dimension $d\geqslant1$,
a vector $\bsk=(k_1,\dots,k_d)$ of integers $k_j\geqslant0$
and a vector $\bsc = (c_1,\dots,c_d)$ of integers with $0\leqslant c_j <b^{k_j}$,
the Cartesian product
$$E(\bsk,\bsc) = \prod_{j=1}^d\Bigl[\frac{c_j}{b^{k_j}}, \frac{c_j+1}{b^{k_j}}\Bigr)
$$
is an elementary interval in base $b$.
It has volume $b^{-|\bsk|}$ where $|\bsk| = \sum_{j=1}^dk_j$.

Speaking informally, the set $E(\bsk,\bsc)$ has a proportion
$b^{-|\bsk|}$ of the volume of $[0,1]^d$ and so it
`deserves' to get (i.e., contain) $nb^{-|\bsk|}$ points when
we place $n$ points inside $[0,1]^d$.
Digital nets satisfy that condition for certain~$\bsk$.
We use the following definitions from Niederreitter \cite{nied:1987}.

\begin{definition}\label{def:tmdnet}
For integers $m\geqslant t\geqslant0$,
the $n=b^m$ points $\bsu_1,\dots,\bsu_n\in[0,1]^d$ are a $(t,m,d)$-net
in base $b\geqslant2$, if every elementary interval $E(\bsk,\bsc)\subset[0,1]^d$
of volume $b^{t-m}$ contains exactly $b^t$ of the points
$\bsu_1,\dots,\bsu_n$.
\end{definition}

Every elementary interval that `deserves' $b^t$ points
of the digital net, gets that many of them.
When we speak of digital nets we ordinarily mean $(t,m,d)$-nets
though some authors reserve the term `digital' to refer to
specific construction algorithms rather than just the property in Definition~\ref{def:tmdnet}.

\begin{definition}\label{def:tdseq}
For integers $t\geqslant0$, $b\geqslant2$ and $d\geqslant1$, the infinite sequence
$\bsu_1,\bsu_2,\dots\in[0,1]^d$ is a $(t,d)$-sequence in base $b$
if $\bsu_{(r-1)b^m+1},\dots,\bsu_{rb^m}$
is a $(t,m,d)$-net in base $b$ for any integers $m\geqslant t$ and $r\geqslant1$.
\end{definition}

Sobol' sequences \cite{sobo:1967:tran} are $(t,d)$-sequences in base $2$.
From Definition~\ref{def:tdseq},
we see that the first $2^m$ points
of a Sobol' sequence are a $(t,m,d)$-net in base $2$
for any $m\geqslant t$.
So are the second $2^m$ points, and if we merge both
of those point sets, we get a $(t,m+1,d)$-net in base $2$.
We can merge the first two of those to get a $(t,m+2,d)$-net
in base $2$ and so on ad infinitum.

Given $b$, $m$ and $d$, smaller values of $t$ are better.
It is not always possible to have $t=0$ and the best
possible $t$ increases monotonically with $d$.
The best known values of $t$ for $(t,d)$-sequences
and $(t,m,d)$-nets are given in the online MinT web site \cite{schu:schm:2009}.
The published $t$ value for a Sobol' sequence might be
conservative in that the first $b^m$ points of the Sobol' sequence
can possibly be a $(t',m,d)$-net for some $t'<t$.

The proven properties of digital nets
including those taken from Sobol' sequences
derive from their balanced sampling of elementary
intervals.
The analysis path can be via discrepancy \cite{nied:1992}
or Haar wavelets \cite{sobo:1969}
or Walsh functions~\cite{dick:pill:2010}.

%Any function that is constant within
%an elementary interval of volume $b^{t-m}$
%will be correctly integrated by averaging over
%a digital net. The same is true for any linear combination
%of such functions taken over the combinatorially many elementary
%intervals that get balanced as $m$ increases.

%\cite{chen:sriv:trav:2014},

\begin{figure}[t!]
\centering
\includegraphics[width=.9\hsize]{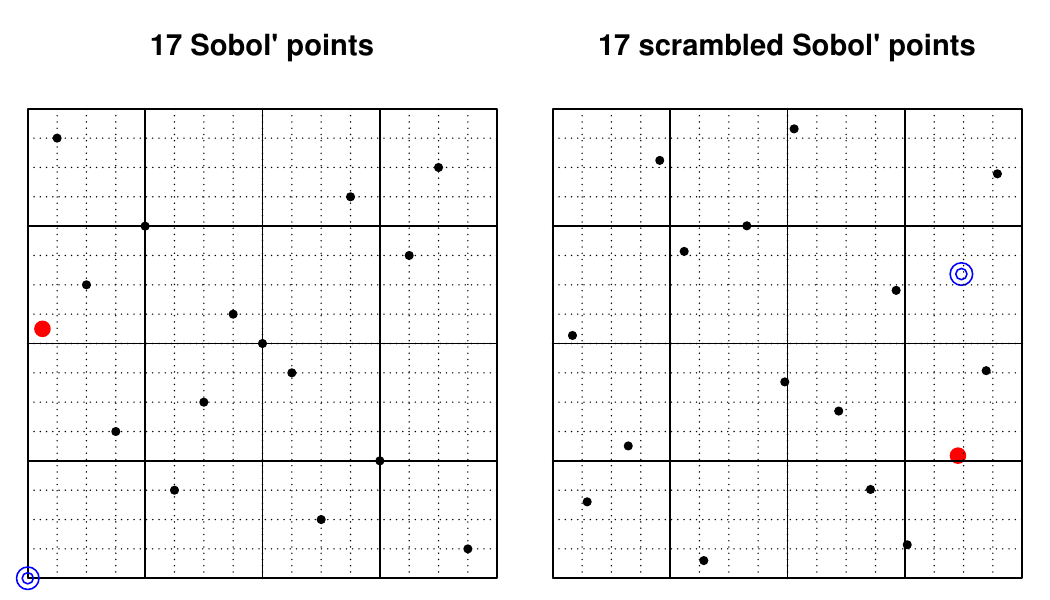}
\caption{\label{fig:notnet}
The left panel shows the first $17$ Sobol' points in $[0,1]^2$.
The initial point at $(0,0)$ is shown in concentric circles.
The $17$'th point is shown as a large disk.
Solid reference lines partition $[0,1]^2$ into $16$
congruent squares. Dashed reference lines partition it
into $256$ congruent squares.  The right panel shows
a nested uniform scramble of these $17$ points.
}
\end{figure}

The left panel in Figure~\ref{fig:notnet}
shows the first $16$ points of a Sobol' sequence
in two dimensions. Fifteen of them are small solid
disks and one other is represented by concentric
circles at the origin.
Those points form a $(0,4,2)$-net in base $2$.
Reference lines divide the unit square into
a $4\times 4$ grid of elementary intervals
of size $1/4\times 1/4$.
Each of those has one of the $16$ points,
often at the lower left corner. Recall that
elementary intervals include their lower boundary
but not their upper boundary.
Finer reference lines partition the unit square into
$16$ strips of size $1\times 1/16$.  Each of those
has exactly one point of the digital net.
The same holds for the $16$ rectangles of
each of these shapes:
$1/2\times 1/8$, $1/8\times 1/2$ and $1/16\times 1$.
All told, those $16$ points have balanced $80$ elementary
intervals and the number of balanced intervals grows
rapidly with $m$ and $d$.

The point $\bsu_1=(0,0)$ is problematic as described
above.  If we skip it and take points $\bsu_2,\dots,\bsu_{17}$
then we replace it with the large solid disk
at $(1/32,17/32)$.  Doing that leaves the lower left
$1/4\times 1/4$ square empty and puts two points
into a square above it.  The resulting $16$ points now
fail to be a $(0,4,2)$-net.

The introduction mentioned some randomizations
of digital nets.  There is a survey of RQMC
in \cite{lecu:lemi:2000}.
For definiteness, we consider
the nested uniform scramble from \cite{rtms}.
Applying a nested uniform scramble to a $(t,d)$-sequence
$\bsu_1,\bsu_2,\dots$ in base $b$ yields points
$\bsx_1,\bsx_2,\dotsi$ that individually satisfy $\bsx_i\sim\dunif[0,1]^d$
and collectively are a $(t,d)$-net in base $b$ with probability one.
The estimate $\hat\mu_{\bsx,1}$ then satisfies
$\e(\hat\mu_{\bsx,1})  = \mu$ by uniformity of $\bsx_i$.
The next paragraphs summarize some additional properties of scrambled nets.

If $f\in L^{1+\epsilon}[0,1]^d$ for some $\epsilon>0$
then \cite{owen:rudo:2020} show that $\Pr(\lim_{m\to\infty}\hat\mu_{\bsx,1}=\mu) =1$,
where the limit is through $(t,m,d)$-nets formed by initial $b^m$
subsequences the $(t,d)$-sequence of $\bsx_i$.
If $f\in L^2[0,1]^d$ then
$\var(\hat\mu_{\bsx,1})=o(1/n)$ as $n=b^m\to\infty$ \cite{snetvar}.
That is, the RMSE is $o(n^{-1/2})$, superior to MC.
Evidence of convergence rates for RQMC better than $n^{-1/2}$ have been seen for some
unbounded integrands from financial problems.
For instance variance reduction factors with respect to MC
have been seen to increase with sample size in \cite{lecu:2009:fin}.

The usual regularity condition for plain MC is that $f(\bsx)$
has finite variance $\sigma^2$ and the resulting RMSE
is $\sigma n^{-1/2}$. When $f\in L^2[0,1]^d$ with variance $\sigma^2$
then  scrambled net sampling with $n=b^m$ satisfies
\begin{align}\label{eq:gammabound}
\rmse(\hat\mu_{\bsx,1})\leqslant \Gamma^{1/2}\sigma n^{-1/2}
\end{align}
for some $\Gamma<\infty$ \cite{snxs}.
For digital nets in base $2$, such as those of Sobol',
it is known that $\Gamma$ is a power
of two no larger than $2^{t+d-1}$ \cite{pan:owen:2021:tr}.
Equation~\eqref{eq:gammabound} describes a worst case
$f\in L^2[0,1]^d$ that maximizes the ratio of the RMSE
for RQMC to that of MC.

The accuracy of QMC points is most commonly described
by a worst case analysis with $|\mu-\hat\mu|=O(n^{-1}\log(n)^{d-1})$
when $f$ has bounded variation in the sense of Hardy and Krause (BVHK).
These powers of $\log(n)$ are not negligible
for practically relevant values of $n$, when $d$ is moderately large.
Then the bound gives a misleadingly pessimistic idea of the accuracy one can expect.
The bound in~\eqref{eq:gammabound} shows that the RMSE
of scrambled nets is at most $\Gamma^{1/2}\sigma/\sqrt{n}$,
a bound with no powers of $\log(n)$.
This holds for $f\in L^2$, which then includes any $f$
in BVHK as well as many others of practical interest, such as some
unbounded integrands.
Note that integrands in BVHK must be bounded and they are also Riemann integrable \cite{owen:rudo:2020}, and so they are in $L^2$.

Under further smoothness conditions on $f$,
$\rmse(\hat\mu_{\bsx,1}) = O(n^{-3/2}(\log n)^{(d-1)/2})$.
This was first noticed in \cite{smoovar} with a correction in \cite{localanti}.
The weakest known sufficient conditions are a generalized Lipschitz condition from \cite{yue:mao:1999}.
The condition in \cite{localanti} is that for any nonempty $u\subseteq\{1,\dots,d\}$
the mixed partial derivative of $f$ taken once with respect to each index $j\in u$
is continuous on $[0,1]^d$.
To reconcile the appearance and non-appearance of logarithmic factors, those
two results give
$\rmse(\hat\mu_{\bsx,1})\leqslant\min( \Gamma^{1/2}\sigma n^{-1/2}, A_n)$
for some sequence $A_n = O(n^{-3/2}\log(n)^{(d-1)/2})$.
The logarithmic factor can degrade the $n^{-3/2}$ rate but only
subject to a cap on performance relative to plain MC.
Finally, Loh \cite{loh:2003} proves a central limit theorem
for $\hat\mu_{\bsx,1}$ when $t=0$.

The right panel of Figure~\ref{fig:notnet} shows a nested uniform
scramble of the points in the left panel.  The problematic point
$\bsu_1$ becomes a uniformly distributed point in the square,
and is no longer on the boundary. If we replace it by $\bsu_{17}$
then just as in the unscrambled case, there is an empty $1/4\times 1/4$
elementary interval, and another one with two points.

There is a disadvantage to $\hat\mu_{\bsx,2}$ compared to $\hat\mu_{\bsx,1}$
when the latter attains a root mean squared error
$O(n^{-3/2+\epsilon})$, for then
\begin{align}\label{eq:spoiledit}
\hat\mu_{\bsx,2} = \hat\mu_{\bsx,1} + \frac1n\big( f(\bsx_{n+1})- f(\bsx_{1})\bigr).
\end{align}
The term
$(f(\bsx_{n+1})- f(\bsx_{1}))/n = O(1/n)$ will ordinarily decay
more slowly than $|\hat\mu_{\bsx,1}-\mu|$.  Then skipping the first
point will actually make the rate of convergence worse.
A similar problem happens if one simply ignores $\bsx_1$
and averages the $n-1$ points $f(\bsx_2)$ through $f(\bsx_n)$.
A related issue is that when equally weighted integration
rules have errors $O(n^{-r})$ for $r>1$, this rate can only
realistically take place at geometrically separated
values of $n$. See~\cite{sobo:1993:b,quadconstraint}.
The higher order digital nets  of \cite{dick:2008}
attain $o(1/n)$ errors under suitable regularity
conditions and their randomizations in \cite{dick:2011}
attain RMSEs of $o(1/n)$.  The argument against
skipping the first point also apply to these methods.

\section{Synthetic examples}\label{sec:synth}

Here we look at some very simple modest dimensional
integrands. They fit into a `best case' case analysis
for integration, motivated as follows.  We
suppose that some sort of function $g(\bsx)$ is
extremely favorable for a method and also that it resembles
the actual integrand. We may write
$$f(\bsx) = g(\bsx) +\err(\bsx).$$
In the favorable cases, $\err$ is small and $g$ is easily
integrated.  For classical quadratures
$g$ may be a polynomial \cite{davrab}.  For digital nets, some functions
$g$ may have rapidly converging Walsh series \cite{dick:pill:2010},
others are sums of functions of only a few variables at a time
\cite{cafl:moro:owen:1997}.
For lattice rules \cite{sloa:joe:1994},
a favorable $g$ has a rapidly converging Fourier series.
The favorable cases work well because
$$
\frac1n\sum_{i=1}^nf(\bsx_i)
=\frac1n\sum_{i=1}^ng(\bsx_i) +\frac1n\sum_{i=1}^n\err(\bsx_i)
$$
with the first term having small error because it is well
suited to the method and the second term having small
error because $\err(\cdot)$ has a small norm and we take
an equal weight sample of it instead of using large weights
of opposite signs.
A good match between method and $g$ saves us the chore
of searching for one or more control variates.
Choosing cases where a method ought to work
is like the positive controls used in experimental science.
We can use them to verify that the method or its
numerical implementation work as expected on the
cases they were designed for.
There can and will be unfavorable cases in practice.
Measuring the sample variance under replication
provides a way to detect that.

Here we consider some cases where scrambled nets should work well.
The first is
\begin{align}\label{eq:f0}
g_0(\bsx) = \sum_{j=1}^d\bigl(e^{x_j}-e+1\bigr),
\end{align}
which clearly has $\mu=0$.
This sum of centered exponentials is smooth and additive.  It is thus very simple for
QMC and RQMC.
It is unlikely that anybody turns to RQMC for this function
but as remarked above the integrand one has may
be close to such a simple function.

\begin{figure}[t!]
% Warning from author to self:  function numbers in the R source code
% have evolved independently of those in the latex source
\centering
\includegraphics[width=.8\hsize]{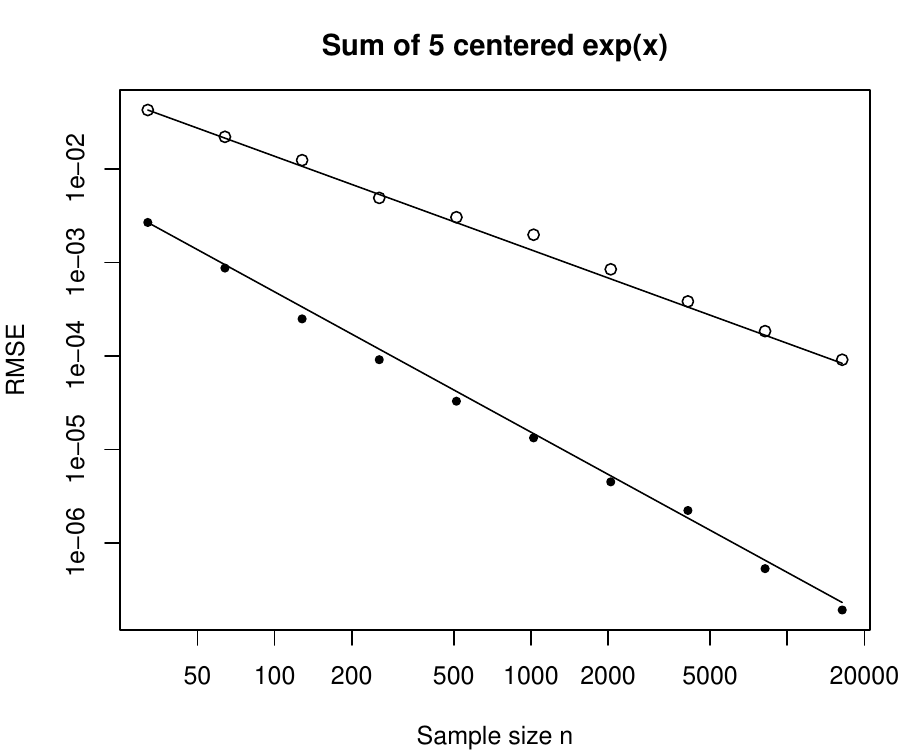}
\caption{\label{fig:f0}
Solid points show RMSE for scrambled Sobol'
estimate $\hat\mu_{\bsx,1}$
versus $n$ from $R=10$ replicates.
A reference line parallel to $n^{-3/2}$ goes through the first solid point.
 Open points show RMSE
for scrambled Sobol' estimates $\hat\mu_{\bsx,2}$ which drop
the initial zero.
A reference line parallel to $n^{-1}$ goes through the first open point.
}
\end{figure}

Figure~\ref{fig:f0} shows the RMSE for
this function $g_0$ based on $R=10$
independent replicates of both $\hat\mu_{\bsx,1}$
and $\hat\mu_{\bsx,2}$.
%Ordinarily $10$ replicates are not enough to estimate a
%mean squared error well but this case is an exception.
%For instance in Figure~\ref{fig:f0} the RMSEs span over 4
%orders of magnitude corresponding to ratios of largest
%to smallest variance greater than $10^8$.
%Even a factor of five in the relative error of estimated
%to true variance would be negligible here.
Reference lines show a clear pattern. The
error follows a reference line parallel to $n^{-3/2}$ on a log-log
plot for $\hat\mu_{\bsx,1}$. For $\hat\mu_{\bsx,2}$, the
reference line is parallel to $n^{-1}$.
These slopes are exactly what we would expect from the underlying theory,
the first from \cite{smoovar} and the second from Equation~\eqref{eq:spoiledit}.
In both cases the line goes through the data for $n=32$ and is
then extrapolated to $n=2^{14}=16{,}384$ with the given slopes.
That is a more severe test for the asymptotic theory than fitting by least
squares would be.  In this instance, the asymptotic theory is
already close to the measurements by $n=32$.
%Both theory and measurements show a strong advantage
%to retaining the scrambled initial point.

An earlier version of this article used $g_0(\bsx)=\sum_{j=1}^dx_j$ instead of
the function $g_0$ above.
The RMSEs for that function also closely follow the predicted rates.
It is not however as good a test case because it is antisymmetric about $\bsx=(1/2,\dots,1/2)$,
meaning that $(g_0(\bsx) + g_0(\tilde\bsx))/2=\mu$ for all $\bsx$, where $\tilde\bsx=1-\bsx$ componentwise.
If we use such an antisymmetric function,
then we will get highly accurate results just from having a nearly antithetic
set of evaluation points that may or may not be equidistributed.

\begin{figure}[t!]
\centering
\includegraphics[width=.8\hsize]{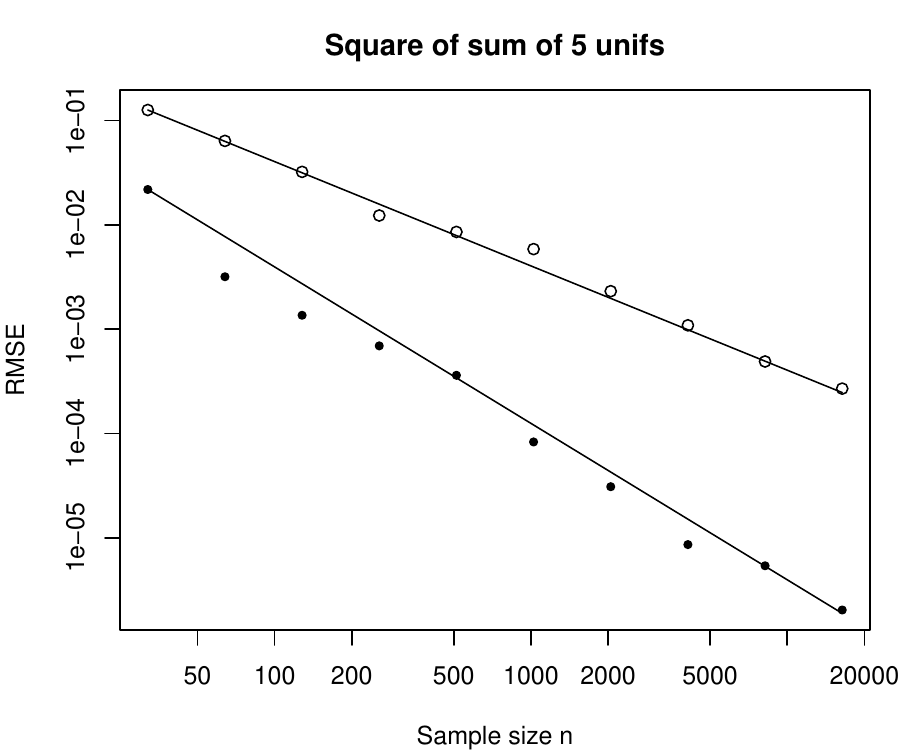}
\caption{\label{fig:f1}
Solid points show RMSE for scrambled Sobol'
estimate $\hat\mu_{\bsx,1}$
versus $n$ from $R=10$ replicates.
A reference line parallel to $n^{-3/2}$ goes through the first solid point.
 Open points show RMSE
for scrambled Sobol' estimates $\hat\mu_{\bsx,2}$ which drop
the initial zero.
A reference line parallel to $n^{-1}$ goes through the first open point.
}
\end{figure}

The second function is
\begin{align}\label{eq:f0}
g_1(\bsx) = \Biggl(\, \sum_{j=1}^dx_j\Biggr)^2.
\end{align}
Unlike $g_0$ this function is not additive.
It has interactions of order $2$ but no higher
in the functional ANOVA decomposition \cite{hoef:1948,sobo:1969}
and it also has a substantial additive component.
It is not antisymmetric about $(1/2,1/2,\dots,1/2)$.
It has $\mu = d/3 + d(d-1)/4$.
Figure~\ref{fig:f1} shows the RMSE for
$\hat\mu_{\bsx,1}$ and  $\hat\mu_{\bsx,2}$.  Once again
they follow reference lines parallel to $n^{-3/2}$ and
$n^{-1}$ respectively.
Asymptotic theory predicts a mean squared error
with a component proportional to $n^{-3}$
and a second component proportional to $\log(n)n^{-3}$
that would eventually dominate the first,
leading to an RMSE that approaches $n^{-3/2}\log(n)^{1/2}$.
%proportional to $\log(n)^{1/2}/n^{3/2}$ arising from the two factor
%interaction but the additive proportion contributes $O(n^{-3/2})$
%with a larger coefficient.

Next we look at a product
$$g_2(\bsx) = \prod_{j=1}^d(e^{x_j}-e+1).$$
This function has $\mu=0$ for any $d$.
It is surprisingly hard for (R)QMC to handle this function
for modest $d$, much less large $d$.  It is dominated by $2^d$ spikes of opposite
signs around the corners of $[0,1]^d$.  It may also be extra hard for Sobol'
points compared to alternatives, because Sobol' points often have
rectangular blocks that alternate between double the uniform density and emptiness.
In a functional ANOVA decomposition,
it is purely $d$-dimensional in that the only non-zero variance
component is the one involving all $d$ variables.
Asymptotic theory predicts an RMSE of $O(n^{-3/2}\log(n)^{(d-1)/2})$.

\begin{figure}[t!]
\centering
\includegraphics[width=.8\hsize]{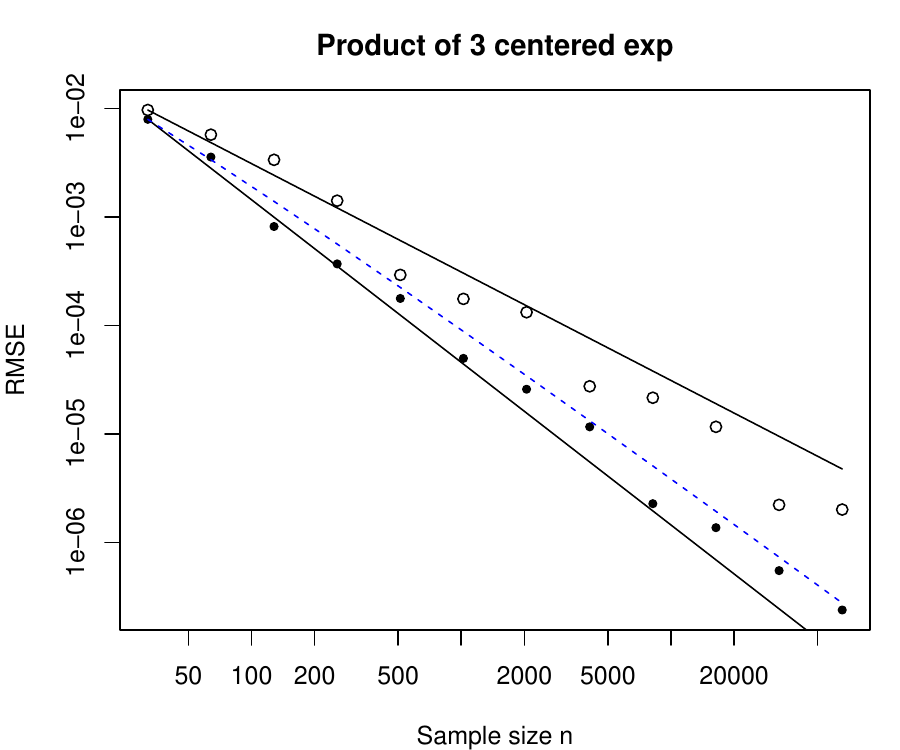}
\caption{\label{fig:f2d3}
The integrand is a product of 3 centered exponentials.
Solid points show RMSE for scrambled Sobol'
estimate $\hat\mu_{\bsx,1}$
versus $n$ from $R=10$ replicates.
A reference line parallel to $n^{-3/2}$ goes through the first solid point.
Open points show RMSE
for scrambled Sobol' estimates $\hat\mu_{\bsx,2}$ which drop
the initial zero.
A reference line parallel to $n^{-1}$ goes through the first open point.
A dashed reference line through the first solid point decays
as $\log(n)/n^{3/2}$.
}
\end{figure}

Figure~\ref{fig:f2d3} shows results for $d=3$ and this $g_2(\bsx)$.
The rate for $\hat\mu_{\bsx,1}$ shows up as slightly
worse than $n^{-3/2}$ while the one for $\hat\mu_{\bsx,2}$
appears to be slighly better than $n^{-1}$.
Both are much better than $O(n^{-1/2})$.
Putting in the predicted logarithmic factor improves the
match between asymptotic prediction and empirical outcome
for $\hat\mu_{\bsx,1}$.
It is not clear what can explain $\hat\mu_{\bsx,2}$ doing
better here than the aysmptotic prediction.
Perhaps the asymptotics become descriptive of actual errors at much larger $n$
for this function than for the others.
Judging by eye it is possible that the convergence rate
is worse when the first point is dropped, but the evidence
is not as clear as in the other figures where the computed
values so closely follow theoretical predictions.
There is an evident benefit to retaining the initial point that at
a minimum manifests as a constant factor of improvement.

In some of the above examples the asymptotic theory
fit very well by $n=32$.
One should not expect this in general.
It is more reasonable to suppose that that is a consequence
of the simple form of the integrands studied in this section.
For these integrands the strong  advantage of retaining
the original point shows in both theory and empirical values.
There is no countervailing theoretical reason to support
dropping the first point.

\section{Wing weight function}\label{sec:wingweight}

The web site \cite{surj:bing:2013}
includes a $10$ dimensional function that
computes the weight of an airplane's wing based on a physical model of the way
the wing is manufactured.  While one does not ordinarily want to know the
average weight of a randomly manufactured wing, this function is interesting
in that it has a real physical world origin instead of being completely
synthetic.  It is easily integrated by several QMC methods
\cite{qmcparts} and so
it is very likely that it equals $g+\err$ for a favorable $g$ and a small~$\err$.

The wing weight function is
\begin{align*}
0.036 S_{\mathrm{w}}^{0.758}W_{\mathrm{fw}}^{0.0035}
\Bigl(\frac{A}{\cos^2(\Lambda)}\Bigr)^{0.6}
q^{0.006}\lambda^{0.04}
\Bigl(\frac{100t_{\mathrm{c}}}{\cos(\Lambda)}\Bigr)^{-0.3}
(N_{\mathrm{x}}W_{\mathrm{dg}})^{0.49}+S_{\mathrm{w}}W_{\mathrm{p}}.
\end{align*}
The definition and uniform ranges of these variables are given in Table~\ref{tab:wwtvariables}.

\begin{table}[t]
\centering
\begin{tabular}{lll}
\toprule
Variable & Range & Meaning\\
\midrule
$S_{\mathrm{w}}$ & [150, 200] & wing area (ft$^2$)\\
$W_{\mathrm{fw}}$ & [220, 300] & weight of fuel in the wing (lb)\\
$A$ & [6, 10] &aspect ratio\\
$\Lambda$ & [$-$10, 10] & quarter-chord sweep (degrees)\\
$q$ & [16, 45] & dynamic pressure at cruise (lb/ft$^2$)\\
$\lambda$ & [0.5, 1] & taper ratio\\
$t_{\mathrm{c}}$ & [0.08, 0.18] & aerofoil thickness to chord ratio\\
$N_{\mathrm{z}}$ & [2.5, 6] &ultimate load factor\\
$W_{\mathrm{dg}}$  & [1700, 2500] & flight design gross weight (lb)\\
$W_{\mathrm{p}}$ & [0.025, 0.08] & paint weight (lb/ft$^2$)\\
\bottomrule
\end{tabular}
\caption{\label{tab:wwtvariables}
Variables and their ranges for the wing weight function.
}
\end{table}

\begin{figure}[t!]
\centering
\includegraphics[width=.8\hsize]{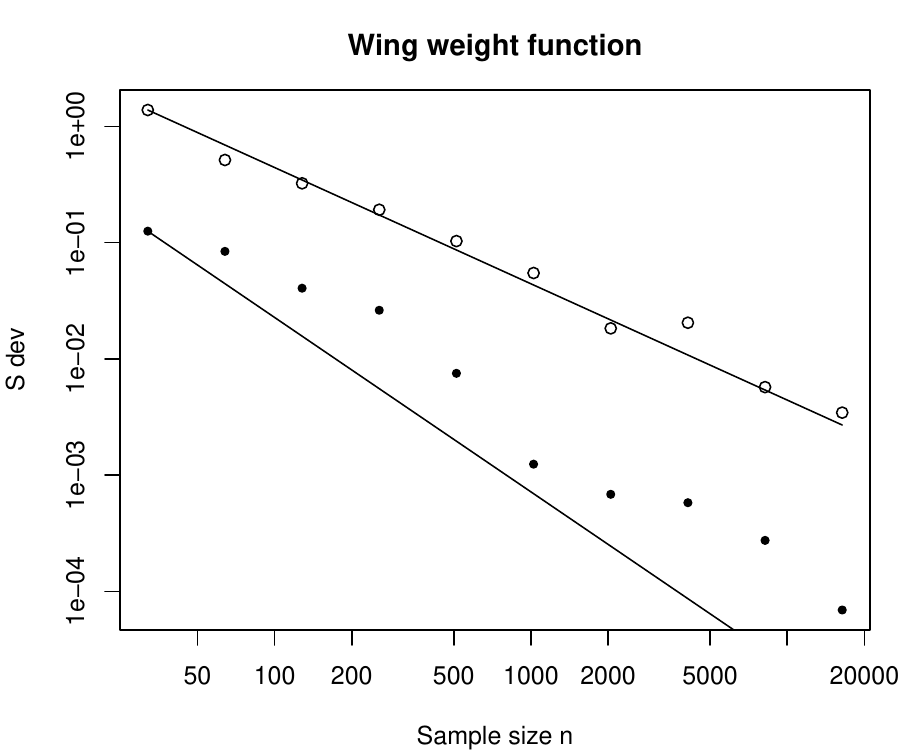}
\caption{\label{fig:outwwt}
Solid points show standard deviation for scrambled Sobol'
estimate $\hat\mu_{\bsx,1}$
versus $n$ from $R=10$ replicates.
A reference line parallel to $n^{-3/2}$ goes through the first solid point.
Open points show standard deviation
for scrambled Sobol' estimates $\hat\mu_{\bsx,2}$ which drop
the initial zero.
A reference line parallel to $n^{-1}$ goes through the first open point.
}
\end{figure}

For this function the standard deviation among $10$ independent replicates
is used instead of the RMSE.
The results are in Figure~\ref{fig:outwwt}.  Once again
there is a strong disadvantage to dropping the first Sobol' point.
The RMSE when dropping the first point is very nearly $O(n^{-1})$.
The RMSE for not dropping the first point is clearly better.
The pattern there is not linear on the log-log scale so
we cannot confidently conclude what convergence
rate best describes it.

\section{Discussion}\label{sec:discussion}

MC and QMC and RQMC points all come as an $n\times d$ matrix of numbers
in $[0,1]$ that we can then pipe through several functions to change
the support set and distribution and finally evaluate a desired integrand.
Despite that similarity, there are sharp difference in the
properties of QMC and RQMC points that affect how we should use them.

This paper has focussed on a small burn-in, dropping just one
of the points and picking up the next $n$. Burn-in makes
no difference to plain MC apart from doing some
unneeded function evaluations, and it can bring large
benefits to MCMC.
See the comment by Neal in the discussion \cite{kass:carl:gelm:neal:1998}.
Burn-in typically spoils the digital net property.  It is safer to scramble the points
which removes the potentially problematic first point at the origin
while also increasing accuracy on very favorable functions
like those in the examples and also on some unfavorable
ones having singularities or other sources of infinite
variation in the sense of Hardy and Krause. See \cite{owen:rudo:2020}.

There are some exceptional cases where burn-in of (R)QMC may
be harmless.  For $d=1$, any consecutive $2^m$ points
of the van der Corput sequence \cite{vand:1935:I} are a $(0,m,1)$-net in base $2$.
As we saw in Figure~\ref{fig:notnet} that is not always true for $d>1$.
Dropping the first $N=2^{m'}$  points
of a Sobol' sequence for $m'\ge m$
should cause no problems because the next $2^m$ points are
still a $(t,m,s)$-net.
Most current implementations of Sobol' sequences are periodic with
$\bsx_i = \bsx_{i+2^M}$ for a value of $M$ that is typically close to $30$.
Then one could take $m'=M-1$ allowing and use
$m$ up to $M-1$.

The Halton sequence \cite{halt:1960} has few if any
especially good sample sizes $n$
and large burn-ins have been used there.
For plain MC points it is natural to use a round number like
$1000$ or $10^6$ of sample points.  That can be very damaging in
(R)QMC if the points were defined for some other sample size.
Using $1000$ points of a Sobol' sequence may well be less accurate than using $512$.
Typical sample sizes are powers of $2$ for digital nets
and large prime numbers for lattice rules \cite{sloa:joe:1994,lecu:lemi:2000}.
The Faure sequences \cite{faur:1982} use $b=p\geqslant d$ where $p$ is
a prime number.
With digital nets as with antibiotics, one should take the whole sequence.

Another practice that works well in MCMC, but should not
be used in (R)QMC is `thinning'.  In MCMC, thinning
can save storage space and in some cases can improve
efficiency despite increasing variance \cite{thinmcmc}.
One takes every $k$'th
point, $\bsx_{k\times i}$ for some integer $k>1$, or in combination with burn-in
$\bsx_{B+k\times i}$ for some integer $B\geqslant1$.
To see the problem, consider the very basic van der Corput
sequence $x_i\in[0,1]$.  If $x_i\in[0,1/2)$ then $x_{i+1}\in[1/2,1)$.
For instance \cite{cafl:mosk:1995} use that observation
to point out that simulating a Markov chain with van der Corput
points can be problematic.
Now suppose that one thins the van der Corput sequence to every second point
using $k=2$.  All of the retained points are then in either $[0,1/2)$
or in $[1/2,1)$. One will estimate either $2\int_0^{1/2}f(x)\rd x$
or $2\int_{1/2}^1f(x)\rd x$ by using that sequence.
The first component of a Sobol' sequence is usually
a van der Corput sequence.

Thinning for QMC was earlier considered by \cite{koci:whit:1997}
who called it `leaping'. They find interesting results taking
every $L$'th point from a Halton sequence, taking $L$ to be
relatively prime to all the bases used in the Halton sequence.
Empirically,
$L=409$ was one of the better values.
They also saw empirically that leaping in digital nets
of Sobol' and Faure lead to non-uniform coverage of the space.

\begin{figure}[t!]
\centering
\includegraphics[width=.8\hsize]{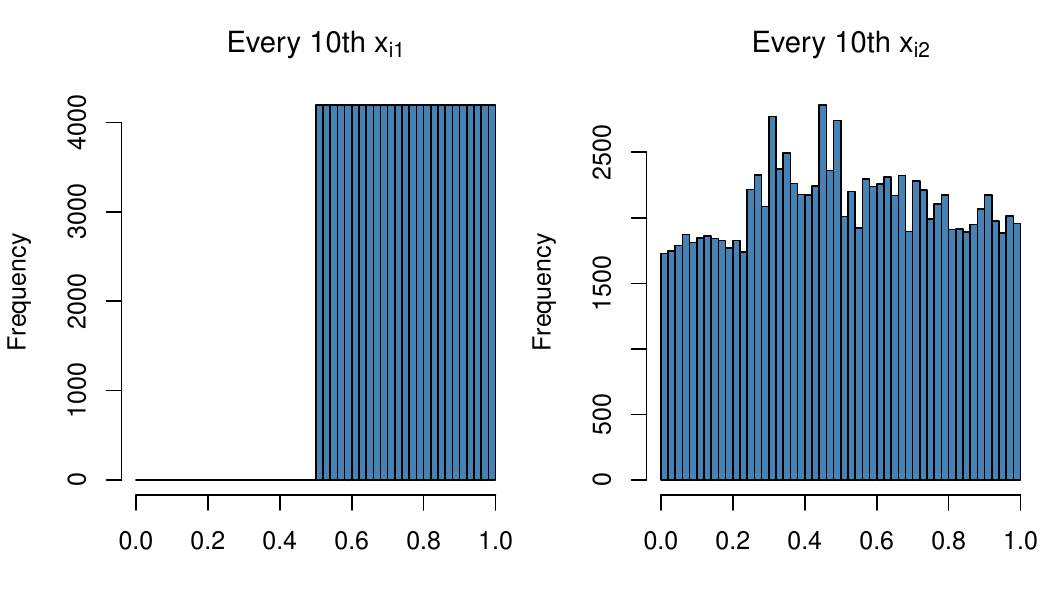}
\caption{\label{fig:thinsobo}
The left panel shows a histogram of every $10$'th
$\bsx_{i1}$ from the first $2^{20}$ Sobol' points.
The right panel shows a histogram of every $10$'th
$\bsx_{i2}$ from the first $2^{20}$ Sobol' points.
}
\end{figure}

The Matlab R2020a {\tt sobolset} function
\url{https://www.mathworks.com/help/stats/sobolset.html}
as of August 11, 2020 includes a thinning/leaping option through
a parameter Leap which is an interval between points,
corresponding to $k-1$ in the discussion above.
It also has a parameter Skip, corresponding to burn-in,
which is a number of initial points to omit.
Fortunately both Leap and Skip are turned off by default.
However even having them present is problematic.
It is not clear how one should use them safely.
The left panel of Figure~\ref{fig:thinsobo} shows a histogram of the values $\bsx_{10i,1}$
for $1\leqslant i\leqslant \lfloor 2^{20}/10\rfloor$.
The right panel  shows a histogram of the values $\bsx_{10i,2}$.

Another area where QMC requires more care than plain MC
is in parallel computing where a task is to be shared over many
processors.
When there are $p$ processors working together, one strategy
from \cite{ParQMC}  is to use a $d+1$ dimensional QMC construction of which
one dimension is used to assign input points to processors.
Processor $k\in\{0,1,\dots,p-1\}$ gets all the points $\bsu_i$ with
$\lfloor px_{i,c}\rfloor=k$ for some $c\in\{1,2,\dots,d+1\}$.
It then uses the remaining $d$ components of $\bsu_i$ in
its computation.  With this strategy each processor gets
a low discrepancy sequence individually which is
better than thinning to every $p$'th point would be.
They collectively have a complete QMC point set.
See \cite{ParQMC} for this and for more references about
parallelizing QMC.

\section*{Acknowledgments}

This work was supported by the NSF under
grant IIS-1837931, and a grant from Hitachi, Ltd. I thank
Fred Hickernell, Pierre L'Ecuyer, Alex Keller, Max Balandat, Michael McCourt,
Pamphile Roy and  Sergei Kucherenko
for stimulating discussions.
I think Sifan Liu for catching an error in some code.
Thanks to Mike Giles, Arnaud Doucet, Alex Keller
and the whole team at ICMS for making MCQMC 2020
happen despite all the pandemic disruption.
This paper benefited from comments of two anonymous
reviewers.  I thank Alex Keller for handling the proceedings volume.
% NSF, Hitachi, Hickernell, Balandat, organizers

\bibliographystyle{spmpsci}
\bibliography{qmc}
\end{document}